\documentclass[11pt]{amsart}         
\DeclareMathOperator{\Hom}{Hom}
\DeclareMathOperator{\Com}{Com}
\DeclareMathOperator{\CCup}{Cup}
\DeclareMathOperator{\I}{I}
\DeclareMathOperator{\1}{id}

\DeclareMathOperator{\dev}{dev}
\DeclareMathOperator{\x}{\times}
\newcommand{\NN}{\mathbb{N}}
\newcommand{\EEnd}{\mathcal End}
\newcommand{\EE}{\mathcal E}
\newcommand{\bul}{\bullet}
\renewcommand{\u}{\smile}
\renewcommand{\=}{:=}
\renewcommand{\t}{\otimes}
\renewcommand{\o}{\circ}

\newcommand{\de}{\delta}
\newcommand{\Ga}{\Gamma}
\newtheorem{thm}{Theorem}[section]
\newtheorem{prop}[thm]{Proposition}

\newtheorem{lem}[thm]{Lemma}
\newtheorem*{Boundary Lemma}{Boundary Lemma}
\newtheorem*{Main Theorem}{Main Theorem}
\newtheorem*{Vertex Proposition}{Vertex Proposition}
\theoremstyle{definition}
 \newtheorem{defn}[thm]{Definition}
\theoremstyle{remark}
 \newtheorem{rem}[thm]{Remark}
\theoremstyle{remark}
 \newtheorem{rems}[thm]{Remarks}
\theoremstyle{definition}
\newtheorem{exam}[thm]{Example}
\theoremstyle{definition}
 \newtheorem{recap}[thm]{Recapitulation}
\theoremstyle{definition}
 \newtheorem{notations}[thm]{Notations}
\numberwithin{equation}{subsection}
\begin{document}
\null
\vskip3true cm

\begin{center}
{\large\bf ON DERIVATION DEVIATIONS \\
          IN AN ABSTRACT PRE-OPERAD}\\
\vskip22true pt
{\large Liivi Kluge and Eugen Paal}
\end{center}

\begin{quote}
\item
{\bf Abstract.} We consider basic algebraic constructions associated with
an \emph{abstract} pre-operad, such as a $\smile$-algebra, total
composition $\bul$, \emph{pre-cobound\-ary} operator $\de$ and
\emph{tribraces} $\{\cdot,\cdot,\cdot\}$. A derivation deviation of the
pre-coboundary operator over the tribraces is calculated in terms of the
$\smile$-multiplication and total composition.

\smallskip
{\bf Classification~(MSC2000).} 18D50.

\end{quote}

\renewcommand{\thefootnote}{}
\footnote{math.QA/0105229}
\footnote{Research supported in part by the ESF grant 3654.}

\section{Introduction and outline of the paper}

\noindent We consider basic algebraic constructions associated with an
\emph{abstract} pre-operad $C$ (Sec.~\ref{pre-operad}), such as a
$\smile$-algebra (Sec.~\ref{cup-mul}), total composition $\bul$
(Sec.~\ref{total}), \emph{pre-coboundary} operator $\de$ (Sec.~\ref{cup and
pre-coboundary}) and \emph{tribraces} $\{\cdot,\cdot,\cdot\}$
(Sec.~\ref{tribraces}). Main result of the present paper is the Main Theorem
(Sec.~\ref{tribraces}). By defining (see degree notations in
Sec.~\ref{pre-operad}) a \emph{derivation deviation}
(Definition~\ref{derdev}) of the pre-coboundary operator $\de$ over the
tribraces $\{\cdot,\cdot,\cdot\}$,
\[
\dev_{\{\cdot,\cdot,\cdot\}}\de\,(h\t f\t g)
    \=\de\{h,f,g\}-\{h,f,\de g\}-(-1)^{|g|}\{h,\de f,g\}
      -(-1)^{|g|+|f|}\{\de h,f,g\},
\]
the Main Theorem states that in a the pre-operad $C$ one has
\[
(-1)^{|g|}\dev_{\{\cdot,\cdot,\cdot\}}\de\,(h\t f\t g)=
    (h\bul f)\u g+(-1)^{|h|f}f\u(h\bul g)-h\bul(f\u g).
\]

This formula (Main Theorem) modifies and generalizes Gerstenhaber's
Theorem~5 from \cite{Ger}.
One can see the appearance of the extra terms with $\de h,\de f$ and $\de g$
inside the tribraces.
The Gerstenhaber formula is obtained when
arguments in tribraces $\{h,f,g\}$ are the Hochshild \emph{cocycles},
i.~e. if $\de h=\de f=\de g=0$.
As it was stated in \cite{GerVor94,VorGer},
on the Hochschild \emph{cochain} level the extra terms with
$\de h,\de f$ and $\de g$ also appear when modifying the Gerstenhaber
formula.
We work out these extra terms for an \emph{abstract} pre-operad
and give their interpretation via a \emph{derivation deviation} of
$\delta$ over the tribraces.

In this paper, we do not assume \emph{associativity} constraint
(see Example~\ref{associativity}).

\section{Pre-operad (composition system)}
\label{pre-operad}

\noindent
Let $K$ be a unital commutative associative ring, and let
$C^n$ ($n\in\NN$) be unital $K$-modules.
For \emph{homogeneous} $f\in C^n$, we refer to $n$ as the \emph{degree}
of $f$ and often write (when it does not cause confusion) $f$ instead of $\deg f$;
for example, $(-1)^f\=(-1)^n$, $C^f\=C^n$ and $\o_f\=\o_n$.
Also, it is convenient to use the shifted (\emph{desuspended}) degree
$|f|\=n-1$.
Throughout this paper, we assume that $\t\=\t_K$.

\begin{defn}
A linear (right) \emph{pre-operad} (\emph{composition system})
with coefficients in $K$ is a sequence $C\=\{C^n\}_{n\in\NN}$ of unital
$K$-modules (an $\NN$-graded $K$-module),
such that the following conditions hold:
\begin{enumerate}
\item
For $0\leq i\leq m-1$ there exist \emph{partial compositions}
\[
  \o_i\in\Hom\,(C^m\t C^n,C^{m+n-1}),\qquad |\o_i|=0.
\]
\item
For all $h\t f\t g\in C^h\t C^f\t C^g$, the \emph{composition relations} hold,
\[
(h\o_i f)\o_j g=
   \begin{cases}
    (-1)^{|f||g|} (h\o_j g)\o_{i+|g|}f
                       &\text{if $0\leq j\leq i-1$},\\
    h\o_i(f\o_{j-i}g)  &\text{if $i\leq j\leq i+|f|$},\\
    (-1)^{|f||g|}(h\o_{j-|f|}g)\o_i f
                       &\text{if $i+f\leq j\leq|h|+|f|$}.
\end{cases}
\]
\item
There exists a unit $\I\in C^1$ such that
\[
\I\o_0 f=f=f\o_i \I,\qquad 0\leq i\leq |f|.
\]
\end{enumerate}
\end{defn}

In the 2nd item, the \emph{first} and \emph{third} parts of the
defining relations turn out to be equivalent.

\begin{rems}
A pre-operad is also called a
              \emph{comp(osition) algebra}
           or \emph{asymmetric operad}
           or \emph{non-symmetric operad}
           or \emph{non-$\Sigma$ operad}.
The concept of (\emph{symmetric}) \emph{operad} was formalized by
May \cite{May72} as a tool for the theory of iterated loop spaces.
Recent studies and applications can be found in \cite{Rene}.

Above we modified the Gerstenhaber \emph{comp algebra} defining relations
\cite{CGS93,GGS92Am} by introducing the sign $(-1)^{|f||g|}$ in the
defining relations of the pre-operad.
The modification enables us to keep track of (control) sign changes more
effectively.
One should also note that (up to sign) our $\o_i$ is Gerstenhaber's $\o_{i+1}$
from \cite{GGS92Am,CGS93}; we use the original (non-shifted)
convention from \cite{Ger,Ger68}.
\end{rems}

\begin{exam}[endomorphism pre-operad {\rm \cite{Ger,Ger68,GGS92Am}}]
\label{HG}
Let $A$ be a unital $K$-module and
$\EE_A^n\={\EEnd}_A^n\=\Hom\,(A^{\t n},A)$.
Define the partial compositions for $f\t g\in\EE_A^f\t\EE_A^g$ as
\[
f\o_i g\=(-1)^{i|g|}f\o(\1_A^{\t i}\t g\t\1_A^{\t(|f|-i)}),
         \qquad 0\leq i\leq |f|.
\]
Then $\EE_A\=\{\EE_A^n\}_{n\in\NN}$ is a pre-operad
(with the unit $\1_A\in\EE_A^1$) called the \emph{endomorphism pre-operad}
of $A$.
A few examples (without the sign factor) can be found in
\cite{Ger68,GGS92Am} as well.
We use the original indexing of \cite{Ger,Ger68} for the
defining formulae.
\end{exam}

\begin{exam}[associahedra]
A geometrical  example of a pre-operad (without signs $(-1)^{...}$ in the
defining relations) is provided by the Stasheff \emph{associahedra}, which
was first constructed in \cite{Sta63}. Quite a surprising realization of
the associahedra as \emph{truncated simplices} was discovered and studied
in \cite{ShnSte94,Sta97,Markl97}.
\end{exam}

\begin{notations}[scope of a pre-operad]
\label{scope}
The {\it scope} of $(h\o_i f)\o_j g$ is given by
\[
0\leq i\leq|h|,\qquad 0\leq j\leq |f|+|h|.
\]
It follows from the defining relations of a pre-operad that
the scope is a disjoint union of
\allowdisplaybreaks
\begin{align*}
B &\=\{(i,j)\in \NN\x\NN\,|\, 1\leq i\leq|h|  \,;\, 0   \leq j\leq i-1\},\\
A &\=\{(i,j)\in \NN\x\NN\,|\, 0\leq i\leq|h|  \,;\, i   \leq j\leq i+|f|\},\\
G &\=\{(i,j)\in \NN\x\NN\,|\, 0\leq i\leq|h|-1\,;\, i+f \leq j\leq |f|+|h|\}.
\end{align*}
\allowdisplaybreaks
Note that the triangles $B$ and $G$ are symmetrically situated
with respect to the parallelogram $A$ in the scope $BAG\=B\sqcup A\sqcup G$.
The (recommended and impressive) picture is left for a reader as
an exercise.
\end{notations}

\begin{recap}
The defining relations of a pre-operad can be easily rewritten as follows:
\[
(h\o_i f)\o_j g=
\begin{cases}
   (-1)^{|f||g|}(h\o_j g)\o_{i+|g|}f &\text{if $(i,j)\in B$},\\
   h\o_i(f\o_{j-i}g)                 &\text{if $(i,j)\in A$},\\
   (-1)^{|f||g|}(h\o_{j-|f|}g)\o_i f &\text{if $(i,j)\in G$}.
\end{cases}
\]
The \emph{first} ($B$) and \emph{third} ($G$) parts of the relations
turn out to be equivalent.
\end{recap}

\section{Cup}
\label{cup-mul}

\noindent
In this section, we recollect from \cite{KPS} basic facts about
cup-multiplication in an \emph{abstract} pre-operad.

\begin{defn}[{\rm cf. \cite{GGS92Am,CGS93}}]
In a pre-operad $C$, let $\mu\in C^2$. Define
$\u\=\u_\mu\:C^f\t C^g\to C^{f+g}$ by
\[
f\u g\=(-1)^f(\mu\o_0 f)\o_f g\in C^{f+g},
\qquad|\smile|=1,\qquad f\t g\in C^f\t C^g.
\]
Note that $\o_f\=\o_{\deg f}$ and $|\smile|=1$ means that $\deg\smile=0$.
The pair $\CCup\,C\=\{C,\u\}$ is called a $\u$-algebra of $C$.
\end{defn}

\begin{exam}
For the endomorphism pre-operad (Example \ref{HG}) $\EE_A$,
one has
\[
f\u g=(-1)^{fg}\mu\o(f\t g),
      \qquad \mu\t f\t g\in \EE_A^2\t \EE_A^f\t \EE_A^g.
\]
\end{exam}

\begin{prop}
\label{cuppro}
In a pre-operad $C$, one has
\[
\mu\o_0 f=(-1)^f f\u\I,
\quad \mu\o_1 f=-\I\u f,
\quad f\u g=-(-1)^{|f|g}(\mu\o_1 g)\o_0f.
\]
\end{prop}
\begin{proof}
We have
\[
(-1)^f f\u\I=(-1)^{f+f}(\mu\o_0 f)\o_f\I=\mu\o_0 f,
           \quad -\I\u f=(\mu\o_0\I)\o_1 f=\mu\o_1 f.
\]
Also, calculate
\allowdisplaybreaks
\begin{align*}
f\u g &=(-1)^f(\mu\o_0 f)\o_f g
       =(-1)^{|f||g|+f}(\mu\o_{f-|f|}g)\o_0 f\\
      &=(-1)^{|f||g|+|f|+1}(\mu\o_1 g)\o_0 f
       =-(-1)^{|f|g}(\mu\o_1 g)\o_0 f,
\end{align*}
\allowdisplaybreaks
which is the required formula.
\end{proof}

\begin{lem}
\label{cup}
In a pre-operad $C$, the following composition relations hold:
\[
(f\u g)\o_j h=
   \begin{cases}
     (-1)^{g|h|}(f\o_j h)\u g &\text{if $0\leq j\leq|f|$},\\
     f\u(g\o_{j-f}h)          &\text{if $f\leq j\leq |g|+f$}.
\end{cases}
\]
\end{lem}
\begin{proof}
Calculate, by using the defining relations of a pre-operad:
\allowdisplaybreaks
\begin{align*}
(f\u g)\o_j h &=(-1)^f[(\mu\o_0 f)\o_f g]\o_j h\\
              &= \begin{cases}
(-1)^{f+|g||h|}[(\mu\o_0 f)\o_j h]\o_{f+|h|}g &\text{if $0\leq j\leq|f|$},\\
(-1)^f(\mu\o_0 f)\o_f (g\o_{j-f}h)            &\text{if $f\leq j\leq |g|+f$},\\
                 \end{cases}\\
              &=\begin{cases}
(-1)^{f+|g||h|}[\mu\o_0(f\o_j h)]\o_{f+|h|}g &\text{if $0\leq j\leq|f|$},\\
f\u(g\o_{j-f}h)                              &\text{if $f\leq j\leq |g|+f$},\\
                \end{cases}\\
              &=\begin{cases}
(-1)^{|f|+|h|+1+f+|g||h|}(f\o_j h)\u g &\text{if $0\leq j\leq|f|$},\\
f\u(g\o_{j-f}h)                        &\text{if $f\leq j\leq |g|+f$},\\
                \end{cases}\\
              &=\begin{cases}
(-1)^{g|h|}(f\o_j h)\u g               &\text{if $0\leq j\leq|f|$},\\
f\u(g\o_{j-f}h)                        &\text{if $f\leq j\leq |g|+f$},
                \end{cases}
\end{align*}
\allowdisplaybreaks
which is the required formula.
\end{proof}

\section{Total composition}
\label{total}

\begin{defn}[{\rm cf. \cite{GGS92Am,CGS93}}]
In a pre-operad $C$, the \emph{total composition}
$\bul\:C^f\t C^g\to C^{f+g-1}$ is defined by
\[
f\bul g\=\sum_{i=0}^{|f|}f\o_i g\in C^{f+g-1},
     \qquad |\bul|=0, \qquad f\t g\in C^f\t C^g.
\]
The pair $\Com\,C\=\{C,\bul\}$ is called the \emph{composition algebra} of $C$.
\end{defn}

\begin{thm}
\label{right der}
In a pre-operad $C$, one has
\[
(f\u g)\bul h=f\u(g\bul h)+(-1)^{|h|g}(f\bul h)\u g.
\]
\end{thm}
\begin{proof}
Use Lemma \ref{cup}. Note that $|f\u g|=f+g-1$ and calculate,
\allowdisplaybreaks
\begin{align*}
(f\u g)\bul h &= \sum_{i=0}^{f+g-1}(f\u g)\o_i h
 =\sum_{i=0}^{f-1}(f\u g)\o_i h+\sum_{i=f}^{f+g-1}(f\u g)\o_i h\\
&=(-1)^{|h|g}\sum_{i=0}^{|f|}(f\o_i h)\u g
  +\sum_{i=f}^{f+g-1}f\u(g\o_{i-f}h)\\
&=(-1)^{|h|g}(f\bul h)\u g+\sum_{i'=0}^{|g|}f\u(g\o_{i'}h)\\
&=(-1)^{|h|g}(f\bul h)\u g+f\u(g\bul h),
\end{align*}
\allowdisplaybreaks
which is the required formula.
\end{proof}

\begin{rem}
This theorem tells us that \emph{right} translations in $\Com\,C$ are
(right) derivations of the $\u$-algebra.
It may be anticipated from Theorem~5 of \cite{Ger} that the \emph{left}
translations in $\Com\,C$ are not derivations of the $\u$-algebra
(see the Main Theorem in Sec. \ref{tribraces}).
\end{rem}

\section{Cup and a pre-coboundary operator}
\label{cup and pre-coboundary}

\begin{defn}
In a pre-operad $C$, define a \emph{pre-coboundary} operator $\de_\mu$ by
\[
-\de_\mu f\=[f,\mu]\=f\bul\mu-(-1)^{|f|}\mu\bul f,
             \qquad \mu\t f\in C^2\t C^f.
\]
\end{defn}

\begin{exam}
\label{associativity}
In the Gerstenhaber theory \cite{Ger}, $C$ is an \emph{endomorphism}
pre-operad and $\de_\mu$ is the Hochschild \emph{coboundary operator} with
the property $\de_\mu^2=0$, the latter is due to the \emph{associativity}
$\mu\bul\mu=0$.

In this paper, we do not assume the associativity constraint.
\end{exam}

\begin{prop}
\label{pre-coboundary}
In a pre-operad $C$, one has
\[
-\de_\mu f=f\u\I+f\bul\mu+(-1)^{|f|}\,\I\u f,
\qquad \mu\t f\in C^2\t C^f.
\]
\end{prop}

\begin{defn}
The \emph{derivation deviation} of $\de_\mu$ over $\bul$ is defined by
\[
\dev_\bul\de_\mu(f\t g)
   \=\de_\mu(f\bul g)-f\bul\de_\mu g-(-1)^{|g|}\de_\mu f\bul g.
\]
\end{defn}

\begin{thm}[{\rm cf. \cite{Ger}}]
\label{second main}
In a pre-operad $C$, one has
\[
(-1)^{|g|}\dev_\bul\de_\mu(f\t g)=f\u g-(-1)^{fg}g\u f,
\qquad \mu\t f\t g\in C^2\t C^f\t C^g.
\]
\end{thm}
\begin{proof}
Two proofs can be found in \cite{KPS}.
\end{proof}

\section{Main theorem and Gerstenhaber's method}
\label{tribraces}

\noindent
In this section, we calculate the derivation deviation of the pre-coboundary
operator over tribraces by using the Gerstenhaber \emph{auxiliary variables method}.
The idea of the method has been illustrated in \cite{KPS}.

\begin{defn}[Gerstenhaber tribraces {\rm (cf. \cite{Ger,CGS93})}]
The Gerstenhaber \emph{tribraces} $\{\cdot,\cdot,\cdot\}$ are
defined as a double sum over the triangle $G$ by
\[
\{h,f,g\}\=\sum_{(i,j)\in G}(h\o_i f)\o_j g\in C^{h+f+g-2},
    \quad |\{\cdot,\cdot,\cdot\}|=0,  \quad h\t f\t g\in C^h\t C^f\t C^g.
\]
\end{defn}

\begin{defn}
\label{derdev}
The \emph{derivation deviation} of $\de\=\de_\mu$ over the tribraces
$\{\cdot,\cdot,\cdot\}$ is defined by
\[
\dev_{\{\cdot,\cdot,\cdot\}}\de\,(h\t f\t g)
    \=\de\{h,f,g\}-\{h,f,\de g\}-(-1)^{|g|}\{h,\de f,g\}
      -(-1)^{|g|+|f|}\{\de h,f,g\}.
\]
\end{defn}

\begin{Main Theorem}[{\rm cf.~\cite{Ger,GerVor94,VorGer}}]
\label{main}
In a pre-operad $C$, one has
\[
(-1)^{|g|}\dev_{\{\cdot,\cdot,\cdot\}}\de\,(h\t f\t g)=
    (h\bul f)\u g+(-1)^{|h|f}f\u(h\bul g)-h\bul(f\u g).
\]
\end{Main Theorem}

\begin{notations}[auxiliary variables {\rm(cf.~\cite{Ger})}]
In a pre-operad $C$, for $(i,j)\in G$ define
\begin{align*}
\Ga_{i+1,j+1}\=
    &-(-1)^{|h|+|f|+|g|}\,\I\u((h\o_i f)\o_j g)\\
    &-(-1)^{|f|+|g|}\sum_{s=0}^{i-1}((h\o_s\mu)\o_{i+1}f)\o_{j+1}g
     +(-1)^{|f|+|g|}(h\o_i(\I\!\u\!f))\o_{j+1}g, \\
\Ga'_{i+1,j+1}\=
    &+(-1)^{|g|}(h\o_i(f\u\I))\o_{j+1}g\\
    &-(-1)^{|f|+|g|}\sum_{s=i+1}^{j-f}((h\o_s\mu)\o_if)\o_{j+1}g
     +(-1)^{|g|}(h\o_i f)\o_j(\I\u g), \\
\Ga''_{i+1,j+1}\=
    &+(h\o_i f)\o_j(g\u\I)\\
    &-(-1)^{|f|+|g|}\sum_{s=j-|f|+1}^{|h|}((h\o_s\mu)\o_if)\o_j g
    -((h\o_if)\o_jg)\u\I.
\end{align*}
\end{notations}

\begin{lem}
\label{first}
In a pre-operad $C$, for $(i,j)\in G$ one has
\[
\de((h\o_if)\o_j g)-(h\o_if)\o_j\de g-(-1)^{|g|}(h\o_i\de f)\o_{j+1}g
=\Ga_{i+1,j+1}+\Ga'_{i+1,j+1}+\Ga''_{i+1,j+1}.
\]
\end{lem}
\begin{proof}
See Appendix A.
\end{proof}

\begin{notations}[truncated envelope of $G'$]
Now, define a (finite) \emph{shifted} Gerstenhaber triangle:
\[
G'\=\{(i,j)\in\NN\x\NN\,|\, 1\leq i\leq|h|\,;\,i+f \leq j\leq f+|h|\}.
\]
Its \emph{envelope} is the triangle
\[
G'_{env}\=\{(i,j)\in\NN\x\NN\,|\, 0\leq i\leq h+1\,;\,i+|f|\leq j\leq f+h\}
=G'\sqcup\partial G'_{env}.
\]
The boundary of the envelope $G'_{env}$ is evidently
$\partial G'_{env}=G'_{env}\setminus G'$.
The \emph{truncated envelope} $G'_{\widetilde{env}}$ of $G'$ is defined by
\emph{removing} its vertices,
\[
G'_{\widetilde{env}}\=G'_{env}\setminus\{(|f|,0);(f+h,0);(h+1,f+h)\}
=G'\sqcup\partial G'_{\widetilde{env}}.
\]
Then, evidently, $G'_{\widetilde{env}}$ as having \emph{six} vertices,
turns out to be a (finite) \emph{hexagon}.
\end{notations}

\begin{lem}
\label{second} In a pre-operad $C$, for
$0\leq i\leq|h|$, $i+f\leq j\leq f+|h|$
one has
\[
(-1)^{|f|+|g|}(\de h\o_if)\o_j g=\Ga_{ij}+\Ga'_{i+1,j}+\Ga''_{i+1,j+1},
\]
by definition for $\Ga_{0j},\Ga'_{i,i+|f|}$, and $\Ga''_{i,f+h}$
(boundary values on $\partial G'_{\widetilde{env}}$).
\end{lem}
\begin{proof}
See Appendix B.
\end{proof}

\begin{Boundary Lemma}
\label{boundary}
In a pre-operad $C$, for $(i,j)\in\partial G'_{\widetilde{env}}$ one has
\begin{alignat*}{6}
&\Ga_{0j}      &&=(-1)^{g+|h|f}f\u(h\o_{j-f}g), &&\qquad f\leq j \leq f+|h|,\\
&\Ga'_{i,i+|f|}&&=(-1)^{|g|}   h\o_{i-1}(f\u g),&&\qquad 1\leq i \leq h,\\
&\Ga''_{i,f+h} &&=(-1)^{g}     (h\o_{i-1}f)\u g,&&\qquad 1\leq i \leq h,
\end{alignat*}
by definition for $\Ga_{0f}$, $\Ga'_{h,h+|f|}$, and $\Ga''_{1,f+h}$
(three vertex values).
\end{Boundary Lemma}
\begin{proof}
See Appendix C.
\end{proof}

\subsection*{Proof of the Main Theorem (Gerstenhaber's method)}

First note that
\[
\{h,\de f,g\}
 =\sum_{i=0}^{|h|-1}\sum_{j=i+\de f}^{|h|+|\de f|}(h\o_i \de f)\o_j g
 =\sum_{i=0}^{|h|-1}\sum_{j=i+f}^{|h|+|f|}(h\o_i \de f)\o_{j+1}g.
\]
By using Lemma~\ref{first} we have
\allowdisplaybreaks
\begin{align*}
\de\{h,f,g\}-\{h,f,\de g\}&-(-1)^{|g|}\{h,\de f,g\}\\
&=\sum_{i=0}^{|h|-1}\sum_{j=i+f}^{|h|+|f|}
         (\Ga_{i+1,j+1}+\Ga'_{i+1,j+1}+\Ga''_{i+1,j+1})\\
&=\sum_{i=1}^{|h|}\sum_{j=i+f}^{f+|h|}(\Ga_{ij}+\Ga'_{ij}+\Ga''_{ij})
 =\sum_{G'}(\Ga_{ij}+\Ga'_{ij}+\Ga''_{ij}).
\end{align*}
\allowdisplaybreaks
Now use Lemma~\ref{second} to see that
\allowdisplaybreaks
\begin{align*}
(-1)^{|f|+|g|}\{\de h,f,g\}
&= \sum_{i=0}^{|\de h|-1}\sum_{j=i+f}^{|\de h|+|f|}(-1)^{|f|+|g|}(\de h\o_i f)\o_j g\\
&= \sum_{i=0}^{|h|}\sum_{j=i+f}^{f+|h|}
            (\Ga_{ij}+\Ga'_{i+1,j}+\Ga''_{i+1,j+1})\\
&= \sum_{i=0}^{|h|}\sum_{j=i+f}^{f+|h|}\Ga_{ij}
        +\sum_{i=1}^{h}\sum_{j=i+|f|}^{f+|h|}\Ga'_{ij}
        +\sum_{i=1}^{h}\sum_{j=i+f}^{f+h}\Ga''_{ij}\\
&= \sum_{i=1}^{|h|}\sum_{j=i+f}^{f+|h|}\Ga_{ij}
  +\sum_{j=f}^{f+|h|}\Ga_{0j}
  +\sum_{i=1}^{|h|}\sum_{j=i+f}^{f+|h|}\Ga'_{ij}
  +\sum_{i=1}^{h}\Ga'_{i,i+|f|}\\
& \quad +\sum_{i=1}^{|h|}\sum_{j=i+f}^{f+|h|}\Ga''_{ij}
  +\sum_{i=1}^{h}\Ga''_{i,f+h}\\
&= \sum_{G'}(\Ga_{ij}+\Ga'_{ij}+\Ga''_{ij})
   \hskip3true pt
    +\,\boxed{ \sum_{j=f}^{f+|h|}\Ga_{0j}
              +\sum_{i=1}^{h}\Ga'_{i,i+|f|}
              +\sum_{i=1}^{h}\Ga''_{i,f+h}}\,.
\end{align*}
\allowdisplaybreaks
One can easily see that the resulting \boxed{\emph{boxed formula}} is a sum
over the \emph{boundary
$\partial G'_{\widetilde{env}}$ of the truncated envelope}
$G'_{\widetilde{env}}$ of $G'$.
By cancelling the sums $\sum_{G'}$ and using the Boundary Lemma,
we finally obtain
\allowdisplaybreaks
\begin{align*}
\dev_{\{\cdot,\cdot,\cdot\}} \de\,(f\t g\t h)
=&-(-1)^{g+|h|f}\sum_{j=f}^{f+|h|} f\u(h\o_{j-f}g)\\
 &-(-1)^{|g|}   \sum_{i=1}^{h} h\o_{i-1}(f\u g)
  -(-1)^{g}     \sum_{i=1}^{h}(h\o_{i-1}f)\u g\\
=-(-1)^{g}&\sum_{i=0}^{|h|}
    \Big[(-1)^{|h|f}f\u(h\o_i g)-h\o_i(f\u g)+(h\o_i f)\u g\Big]\\
=(-1)^{|g|}\Big[&(-1)^{|h|f}f\u(h\bul g)-h\bul(f\u g)+(h\bul f)\u g\Big],
\end{align*}
\allowdisplaybreaks
which is the required formula. \qed

\begin{thm}[{\rm cf.~\cite{Ger,VorGer}}]
In a pre-operad $C$, one has
\[
(-1)^{|g|}\dev_{\{\cdot,\cdot,\cdot\}}\de\,(h\t f\t g)=
    [h,f]\u g+(-1)^{|h|f}f\u[h,g]-[h,f\u g].
\]
\end{thm}
\begin{proof}
Combine the Main Theorem with Theorem \ref{right der}.
\end{proof}

\begin{rems}
A well known first form of this theorem, found by Gerstenhaber \cite{Ger}
for the Hochschild \emph{cocycles}, can be seen as a
starting point of the modern
\emph{mechanical mathematics} based nowadays on the pioneering concept of
(\emph{homotopy}) \emph{Gerstenhaber algebra} \cite{GerVor94,VorGer,Vor99}.

Our Main Theorem generalizes and modifies Gerstenhaber's
Theorem~5 from \cite{Ger}.
One can see the appearance of the extra terms with $\de h,\de f$ and $\de g$
inside the tribraces.
The Gerstenhaber formula is obtained when
arguments in tribraces $\{h,f,g\}$ are the Hochshild \emph{cocycles},
i.~e. if $\de h=\de f=\de g=0$.
As it was stated in \cite{GerVor94,VorGer},
on the Hochschild \emph{cochain} level the extra terms
with $\de f,\de g$ and $\de h$ also appear when modifying the
Gerstenhaber formula.
Our contribution is to work out these extra terms for an \emph{abstract}
pre-operad and give their interpretation via a \emph{derivation deviation}
of $\delta$ over the tribraces.
\end{rems}

\section{Appendix A}
\subsection*{Proof of Lemma~\ref{first}}

First note that
%
\begin{align*}
-\de_\mu((h\o_i f)\o_j g)=
           &(-1)^{|h|+|f|+|g|}\,\I\u((h\o_i f)\o_j g)\\
           &+\sum_{s=0}^{|h|+|f|+|g|}((h\o_i f)\o_j g)\o_s\mu
            +((h\o_i f)\o_j g)\u\I.
\end{align*}
%
By using the composition relations
\[
((h\o_i f)\o_j g)\o_s\mu=
\begin{cases}
 (-1)^{|g|}((h\o_i f)\o_s\mu)\o_{j+1}g
                              &\text{if $0\leq s\leq j-1$},  \\
 (h\o_i f)\o_j(g\o_{s-j}\mu)
                              &\text{if $j\leq s\leq j+|g|$},\\
 (-1)^{|g|}((h\o_i f)\o_{s-|g|}\mu)\o_j g
                              &\text{if $j+g\leq s\leq |h|+|f|+|g|$},\\
\end{cases}
\]
cut the above sum in three pieces,
\allowdisplaybreaks
\begin{align*}
-\de((h\o_i f)\o_j g)
  =\,&(-1)^{|h|+|f|+|g|}\,\I\u((h\o_i f)\o_j g)\\
     &+(-1)^{|g|}\sum_{s=0}^{j-1}((h\o_i f)\o_s\mu)\o_{j+1}g
      +\sum_{s=j}^{j+|g|}(h\o_i f)\o_j(g\o_{s-j}\mu)\\
     &+(-1)^{|g|}\sum_{s=j+g}^{|h|+|f|+|g|}((h\o_i f)\o_{s-|g|}\mu)\o_j g
      +((h\o_i f)\o_j g)\u\I.
\end{align*}
\allowdisplaybreaks
Next, cut the first sum $\sum_{s=0}^{j-1}$ in three pieces once again
by using the composition relations
\[
(h\o_i f)\o_s\mu=
\begin{cases}
 (-1)^{|f|}(h\o_s\mu)\o_{i+1}f    &\text{if $0\leq s\leq i-1$},  \\
 h\o_i(f\o_{s-i}\mu)              &\text{if $i\leq s\leq i+|f|$},\\
 (-1)^{|f|}(h\o_{s-|f|}\mu)\o_i f &\text{if $i+f\leq s\leq j-1$},
\end{cases}
\]
and desuspend summation ranges,
\allowdisplaybreaks
\begin{gather*}
\sum_{s=j}^{j+|x|} x\o_{s-j}\mu
  =\sum_{s=0}^{|x|}x\o_s\mu
  =x\bul\mu
  =-(-1)^{|x|}\,\I\u x-\de_\mu x-x\u\I,\quad \text{for}\,\, x=f, g, \\
\sum_{s=i+f}^{j-1} h\o_{s-|f|}\mu
 =\sum_{s=i+1}^{j-f} h\o_{s}\mu,\\
\allowdisplaybreaks
\begin{align*}
\sum_{s=j+g}^{|h|+|f|+|g|}(h\o_i f)\o_{s-|g|}\mu
  &=\sum_{s=j+1}^{|h|+|f|}(h\o_i f)\o_{s}\mu
   =(-1)^{|f|}\sum_{s=j+1}^{|h|+|f|}(h\o_{s-|f|}\mu)\o_i f\\
  &=(-1)^{|f|}\sum_{s=j-|f|+1}^{|h|}(h\o_{s}\mu)\o_i f,
\end{align*}
\allowdisplaybreaks
\end{gather*}
\allowdisplaybreaks
to obtain the required formula. \qed

\section{Appendix B}

\subsection{Proof of Lemma~\ref{second}}
\label{comid}

First note that
\allowdisplaybreaks
\begin{align*}
\Ga_{ij}+\Ga'_{i+1,j}&+\Ga''_{i+1,j+1}
    =-(-1)^{|h|+|f|+|g|}\,\I\u((h\o_{i-1}f)\o_{j-1}g)\\
    &-(-1)^{|f|+|g|}\sum_{s=0}^{i-2}((h\o_s\mu)\o_{i}f)\o_{j}g
     +(-1)^{|f|+|g|}(h\o_{i-1}(\I\u f))\o_{j}g \\
    &+(-1)^{|g|}(h\o_i(f\u\I))\o_{j}g
     -(-1)^{|f|+|g|}\sum_{s=i+1}^{j-f-1}((h\o_s\mu)\o_{i}f)\o_{j}g\\
    &+(-1)^{|g|}(h\o_i f)\o_{j-1}(\I\u g)
     +(h\o_i f)\o_j(g\u\I)\\
    &-(-1)^{|f|+|g|}\sum_{s=j-|f|+1}^{|h|}((h\o_s\mu)\o_{i}f)\o_j g
     -((h\o_if)\o_jg)\u\I.
\end{align*}
\allowdisplaybreaks
We must compare it term by term with
\allowdisplaybreaks
\begin{align*}
-(\de h&\o_i f)\o_j g
 =\,\bigg(
      \Big((-1)^{|h|}\,\I\u h+\sum_{s=0}^{|h|}h\o_s\mu+h\u\I\Big)\o_i f
   \bigg)\o_j g\\
 =&\,(-1)^{|h|}((\I\u h)\o_i f)\o_j g
   +\sum_{s=0}^{i-2}((h\o_s\mu)\o_i f)\o_j g
   +((h\o_{i-1}\mu)\o_i f)\o_j g \\
  &+((h\o_{i}\mu)\o_i f)\o_j g
   +\sum_{s=i+1}^{j-f-1}((h\o_s\mu)\o_i f)\o_j g
   +((h\o_{j-f}\mu)\o_i f)\o_j g \\
  &+((h\o_{j-|f|}\mu)\o_i f)\o_j g
   +\sum_{s=j-|f|+1}^{|h|}((h\o_s\mu)\o_i f)\o_j g
   +((h\u\I)\o_i f)\o_j g.
\end{align*}
\allowdisplaybreaks
Now, recall the sign $(-1)^{|f|+|g|}$ and use composition relations to
note the {\bf\emph{ground identities}}
\allowdisplaybreaks
\begin{gather*}
 ((\I\u h)\o_i f)\o_j g=((\I\u(h\o_{i-1}f))\o_j g
                       =\I\u((h\o_{i-1}f)\o_{j-1}g),\\
 (h\o_{i-1}\mu)\o_i f=h\o_{i-1}(\mu\o_1 f)=-h\o_{i-1}(\I\u f), \\
 (h\o_{i}\mu)\o_i f=h\o_{i}(\mu\o_0 f)=(-1)^{f}h\o_{i}(f\u\I),\\
\begin{align*}
 ((h\o_{j-f}\mu)\o_i f)\o_j g
      =& (-1)^{|f|}((h\o_{i}f)\o_{j-1}\mu)\o_j g
      =  (-1)^{|f|}(h\o_{i}f)\o_{j-1}(\mu\o_1 g)\\
      =& -(-1)^{|f|}(h\o_{i}f)\o_{j-1}(\I\u g),
\end{align*} \\
\begin{align*}
((h\o_{j-|f|}\mu)\o_i f)\o_j g
      =& (-1)^{|f|}((h\o_{i}f)\o_j\mu)\o_j g
      =  (-1)^{|f|}(h\o_{i}f) \o_j(\mu\o_0 g) \\
      =& -(-1)^{|f|+|g|}(h\o_{i}f)\o_j(g\u\I),
\end{align*} \\
 ((h\u\I)\o_i f)\o_j g=(-1)^{|f|}((h\o_i f)\u\I)\o_j g
                      =(-1)^{|f|+|g|}((h\o_i f)\o_j g)\u\I,
\end{gather*}
\allowdisplaybreaks
which prove the required formula.
\qed

\subsection{Proposition/recapitulation}
\emph
{
In a pre-operad $C$, for $(i,j)\in G'$ the auxiliary variables read
}
\begin{align*}
\Ga_{ij}=
    &-(-1)^{|h|+|f|+|g|}((\I\u h)\o_i f)\o_j g
     -(-1)^{|f|+|g|}\sum_{s=0}^{i-1}((h\o_s\mu)\o_{i}f)\o_{j}g,\\
\Ga'_{ij}=
    &-(-1)^{|f|+|g|}\sum_{s=i-1}^{j-f}((h\o_s\mu)\o_{i-1}f)\o_{j}g,\\
\Ga''_{ij}=
    &-(-1)^{|f|+|g|}\sum_{s=j-f}^{|h|}((h\o_s\mu)\o_{i-1}f)\o_{j-1}g
     -(-1)^{|f|+|g|}((h\u\I)\o_{i-1}f)\o_{j-1}g.
\end{align*}
\begin{proof}
Use the {\bf\emph{ground identities}} from the previous subsection
\ref{comid}.
\end{proof}

\section{Appendix C}
\subsection*{Proof of the Boundary Lemma}

First prove the
\begin{Vertex Proposition}
\emph{In a pre-operad $C$, one has}
\begin{alignat*}{2}
 &\Ga_{0,f+|h|}      &&=(-1)^{g+|h|f}f\u(h\o_{|h|}g), \\
 &\Ga'_{1f}         &&=(-1)^{|g|}   h\o_{0}(f\u g),\\
 &\Ga''_{h,f+h}      &&=(-1)^{g}     (h\o_{|h|}f)\u g.
\end{alignat*}
\end{Vertex Proposition}
\begin{proof}
We calculate these vertex values in a standard way, by using Lemma
\ref{second}.
First calculate $\Ga_{0,f+|h|}$.
Use Lemma \ref{second} and $\Ga''_{1,f+h}$  from the Boundary Lemma
to note that
\begin{align*}
\Ga_{0,f+|h|}+\Ga'_{1,f+|h|}+\Ga''_{1,f+h}
  =\Ga_{0,f+|h|}
    &-(-1)^{|f|+|g|}\sum_{s=0}^{|h|}((h\o_s\mu)\o_{0}f)\o_{f+|h|}g\\
    &+(-1)^{g}(h\o_0 f)\u g.
\end{align*}
We must compare it term by term with
\allowdisplaybreaks
\begin{align*}
-(\de h\o_0 f)\o_{f+|h|}g
 &=\,\bigg(
      \Big((-1)^{|h|}\,\I\u h+\sum_{s=0}^{|h|}h\o_s\mu+h\u\I\Big)\o_0 f
    \bigg)\o_{f+|h|}g\\
 &=(-1)^{|h|}((\I\u h)\o_0 f)\o_{f+|h|}g
  +\sum_{s=0}^{|h|}((h\o_s\mu)\o_0 f)\o_{f+|h|}g\\
 &\qquad+((h\u\I)\o_0 f)\o_{f+|h|}g.
\end{align*}
\allowdisplaybreaks
Now, recall the sign $(-1)^{|f|+|g|}$ and use composition relations to note
that
\allowdisplaybreaks
\begin{gather*}
\begin{align*}
(-1)^{|h|}((\I\u h)\o_0 f)\o_{f+|h|}g
   &=(-1)^{|h|+|f|h}((\I\o_0 f)\u h)\o_{f+|h|}g\\
    =(-1)^{|h|f+|f|} (f\u h)\o_{f+|h|}g
   &=-(-1)^{|g|+|f|}\, \boxed{(-1)^{g+|h|f}f\u(h\o_{|h|}g)} \, ,
\end{align*}\\
\begin{align*}
((h\u\I)\o_0 f)\o_{f+|h|}g
  &=(-1)^{|f|}((h\o_0 f)\u\I)\o_{f+|h|}g
   =(-1)^{|f|}(h\o_0 f)\u(\I\o_{0}g)\\
  &=(-1)^{|f|}(h\o_0 f)\u g,
\end{align*}
\end{gather*}
\allowdisplaybreaks
which lead one to the required formula for $\Ga_{0,f+|h|}$.

Next calculate $\Ga'_{1f}$.
Use Lemma \ref{second} and $\Ga_{0,f}$ from the Boundary Lemma to note that
\begin{align*}
\Ga_{0f}+\Ga'_{1f}&+\Ga''_{1,f+1}
            =(-1)^{f|h|+g}f\u(h\o_{0}g)
             +\Ga'_{1f}\\
     &-(-1)^{|f|+|g|}\sum_{s=1}^{|h|}((h\o_s\mu)\o_{0}f)\o_{f}g
      -(-1)^{|f|+|g|}((h\u\I)\o_0 f)\o_f g.
\end{align*}
We must compare it term by term with
\allowdisplaybreaks
\begin{align*}
-(\de h\o_0 f)\o_{f}g
 &=\bigg(
      \Big((-1)^{|h|}\I\u h+\sum_{s=0}^{|h|}h\o_s\mu+h\u\I\Big)\o_0 f
     \bigg)\o_{f}g\\
&=(-1)^{|h|}((\I\u h)\o_0 f)\o_{f}g
  +((h\o_0\mu)\o_0 f)\o_{f}g\\
&\quad+\sum_{s=1}^{|h|}((h\o_s\mu)\o_0 f)\o_{f}g
      +((h\u\I)\o_0 f)\o_{f}g.
\end{align*}
\allowdisplaybreaks
Now, recall the sign $(-1)^{|f|+|g|}$ and use composition relations to note
that
\allowdisplaybreaks
\begin{gather*}
\begin{align*}
((\I\u h)\o_0 f)\o_{f}g
   &=(-1)^{|f|h}((\I\o_0 f)\u h)\o_{f}g
    =(-1)^{|f|h}(f\u h)\o_{f}g\\
   &=(-1)^{|f|h}f\u(h\o_{0}g),
\end{align*}\\
\begin{align*}
((h\o_0\mu)\o_0 f)\o_f g
&=(h\o_0(\mu\o_0 f))\o_f g
=h\o_0((\mu\o_0 f)\o_f g)\\
&=-(-1)^{|f|+|g|}\, \boxed{(-1)^{|g|}h\o_0(f\u g)}\, ,
\end{align*}
\end{gather*}
\allowdisplaybreaks
which lead one to the required formula for $\Ga'_{1f}$.

At last calculate $\Ga''_{h,f+h}$.
Use Lemma \ref{second} and $\Ga'_{h,f+|h|}$  from the Boundary Lemma
to note that
\begin{align*}
\Ga_{|h|,f+|h|}&+\Ga'_{h,f+|h|}+\Ga''_{h,f+h}
=-(-1)^{|f|+|g|+|h|}((\I\u h)\o_{|h|}f)\o_{f+|h|}g\\
&-(-1)^{|f|+|g|}\sum_{s=0}^{|h|-1}((h\o_s\mu)\o_{|h|}f)\o_{f+|h|}g
+(-1)^{|g|}h\o_{|h|}(f\u g)
+\Ga''_{h,f+h}.
\end{align*}
We must compare it term by term with
\allowdisplaybreaks
\begin{align*}
-(\de h\o_{|h|}f)&\o_{f+|h|}g
 =\bigg(
      \Big((-1)^{|h|}\I\u h+\sum_{s=0}^{|h|}h\o_s\mu+h\u\I\Big)\o_{|h|}f
     \bigg)\o_{f+|h|}g\\
 =&\,(-1)^{|h|}((\I\u h)\o_{|h|}f)\o_{f+|h|}g
  +\sum_{s=0}^{|h|-1}((h\o_s\mu)\o_{|h|}f)\o_{f+|h|}g\\
   &+((h\o_{|h|}\mu)\o_{|h|}f)\o_{f+|h|}g
    +((h\u\I)\o_{|h|}f)\o_{f+|h|}g.
\end{align*}
\allowdisplaybreaks
Now, recall the sign $(-1)^{|f|+|g|}$ and use composition relations to note
that
\allowdisplaybreaks
\begin{gather*}
\begin{align*}
((h\o_{|h|}\mu)\o_{|h|}f)\o_{f+|h|}g
&=(h\o_{|h|}(\mu\o_{0}f))\o_{f+|h|}g
=h\o_{|h|}((\mu\o_0 f)\o_f g)\\
&=(-1)^{f}h\o_{|h|}(f\u g),
\end{align*}\\
\begin{align*}
((h\u\I)\o_{|h|}f)&\o_{f+|h|}g
=(-1)^{|f|}((h\o_{|h|}f)\u\I)\o_{f+|h|}g \\
&=(-1)^{|f|}(h\o_{|h|}f)\u(\I\o_{0}g)
=-(-1)^{|f|+|g|}\,\boxed{(-1)^{g}(h\o_{|h|}f)\u g}\, ,
\end{align*}
\end{gather*}
\allowdisplaybreaks
which lead one to the required formula for $\Ga''_{h,f+h}$.
\end{proof}

Now, when the vertex values found, we can accomplish the proof of
the Boundary Lemma.
First calculate $\Ga_{0j}$. Use Lemma \ref{second} to note that
\allowdisplaybreaks
\begin{align*}
\Ga_{0j}+\Ga'_{1j}&+\Ga''_{1,j+1}
    =\Ga_{0j}
      -(-1)^{|f|+|g|}\sum_{s=0}^{j-f}((h\o_s\mu)\o_{0}f)\o_{j}g\\
      &-(-1)^{|f|+|g|}\sum_{s=j-|f|}^{|h|}((h\o_s\mu)\o_{0}f)\o_j g
      -(-1)^{|f|+|g|}((h\u\I)\o_0 f)\o_jg.
\end{align*}
\allowdisplaybreaks
We must compare it term by term with
\allowdisplaybreaks
\begin{align*}
-(\de h\o_0 f)\o_j g
 &=\bigg(
      \Big((-1)^{|h|}\,\I\u h+\sum_{s=0}^{|h|}h\o_s\mu+h\u\I\Big)\o_0 f
   \bigg)\o_j g\\
 &=\,(-1)^{|h|}((\I\u h)\o_0 f)\o_j g
  +\sum_{s=0}^{j-f}((h\o_s\mu)\o_0 f)\o_j g \\
  &\quad+\sum_{s=j-|f|}^{|h|}((h\o_s\mu)\o_0 f)\o_j g
     +((h\u\I)\o_0 f)\o_j g.
\end{align*}
\allowdisplaybreaks
Now, recall the sign $(-1)^{|f|+|g|}$ and use composition relations to
note that
\begin{gather*}
 (-1)^{|h|}((\I\u h)\o_0 f)\o_j g
   =(-1)^{|h|+h|f|}((\I\o_0 f)\u h)\o_{j}g\\
   =(-1)^{|h|f+|f|}(f\u h)\o_{j}g
   =-(-1)^{|f|+|g|}\, \boxed{(-1)^{g+|h|f}f\u(h\o_{j-f}g)}\, ,
\end{gather*}
which leads one to the required formula for $\Ga_{0j}$.

Next calculate $\Ga'_{i,i+|f|}$. Use Lemma \ref{second} to note that
\allowdisplaybreaks
\begin{align*}
\Ga_{i-1,i+|f|}&+\Ga'_{i,i+|f|}+\Ga''_{i,i+f}
 =-(-1)^{|f|+|g|+|h|}((\I\u h)\o_{i-1}f)\o_{i+|f|}g \\
  &-(-1)^{|f|+|g|}\sum_{s=0}^{i-2}((h\o_s\mu)\o_{i-1}f)\o_{i+|f|}g
   +\Ga'_{i,i+|f|} \\
  & -(-1)^{|f|+|g|}\sum_{s=i}^{|h|}((h\o_s\mu)\o_{i-1}f)\o_{i+|f|}g\\
  &-(-1)^{|f|+|g|}((h\u\I)\o_{i-1}f)\o_{i+|f|}g.
\end{align*}
\allowdisplaybreaks
We must compare it term by term with
\allowdisplaybreaks
\begin{align*}
-(\de h\o_{i-1}f)\o_{i+|f|}g
 &=\bigg(
      \Big((-1)^{|h|}\,\I\u h+\sum_{s=0}^{|h|}h\o_s\mu+h\u\I\Big)\o_{i-1}f
   \bigg)\o_{i+|f|}g\\
 &=(-1)^{|h|}((\I\u h)\o_{i-1}f)\o_{i+|f|}g
   +\sum_{s=0}^{i-2}((h\o_s\mu)\o_{i-1}f)\o_{i+|f|}g\\
  &\quad+((h\o_{i-1}\mu)\o_{i-1}f)\o_{i+|f|}g
   +\sum_{s=i}^{|h|}((h\o_s\mu)\o_{i-1}f)\o_{i+|f|}g \\
  &\quad+((h\u\I)\o_{i-1}f)\o_{i+|f|}g.
\end{align*}
\allowdisplaybreaks
Now, recall the sign $(-1)^{|f|+|g|}$ and use composition relations to
note that
\allowdisplaybreaks
\begin{align*}
((h\o_{i-1}\mu)\o_{i-1}f)\o_{i+|f|}g
&=(h\o_{i-1}(\mu\o_{0}f))\o_{i+|f|}g
 =h\o_{i-1}((\mu\o_{0}f)\o_{f}g)\\
&=-(-1)^{|f|+|g|}\,\boxed{(-1)^{|g|}h\o_{i-1}(f\u g)}\,,
\end{align*}
\allowdisplaybreaks
which leads one to the required formula for $\Ga'_{i,i+|f|}$.

At last calculate $\Ga''_{i,f+h}$. Use Lemma \ref{second} to note that
\allowdisplaybreaks
\begin{align*}
\Ga_{i-1,f+|h|}&+\Ga'_{i,f+|h|}+\Ga''_{i,f+h}
 =-(-1)^{|f|+|g|+|h|}((\I\u h)\o_{i-1}f)\o_{f+|h|}g \\
  &-(-1)^{|f|+|g|}\sum_{s=0}^{i-2}((h\o_s\mu)\o_{i-1}f)\o_{f+|h|}g \\
  &-(-1)^{|f|+|g|}\sum_{s=i-1}^{|h|}((h\o_s\mu)\o_{i-1}f)\o_{f+|h|}g
   +\Ga''_{i,f+h}.
\end{align*}
\allowdisplaybreaks
We must compare it term by term with
\allowdisplaybreaks
\begin{align*}
-(\de h\o_{i-1}f)&\o_{f+|h|}g
 =\bigg(
      \Big((-1)^{|h|}\,\I\u h+\sum_{s=0}^{|h|}h\o_s\mu+h\u\I\Big)\o_{i-1}f
   \bigg)\o_{f+|h|}g\\
 &=(-1)^{|h|}((\I\u h)\o_{i-1}f)\o_{f+|h|}g
   +\sum_{s=0}^{i-2}((h\o_s\mu)\o_{i-1}f)\o_{f+|h|}g\\
  &\quad+\sum_{s=i-1}^{|h|}((h\o_s\mu)\o_{i-1}f)\o_{f+|h|}g
  +((h\u\I)\o_{i-1}f)\o_{f+|h|}g.
\end{align*}
\allowdisplaybreaks
Now, recall the sign $(-1)^{|f|+|g|}$ and use composition relations to
note that
\allowdisplaybreaks
\begin{align*}
((h\u\I)&\o_{i-1}f)\o_{f+|h|}g
=(-1)^{|f|}((h\o_{i-1}f)\u\I)\o_{f+|h|}g\\
&=(-1)^{|f|}(h\o_{i-1}f)\u(\I\o_{0}g)
=-(-1)^{|f|+|g|}\,\boxed{(-1)^{g}(h\o_{i-1}f)\u g}\, ,
\end{align*}
\allowdisplaybreaks
which leads one to the required formula for  $\Ga''_{i,f+h}$.
\qed

\section*{Acknowledgements}
\noindent
One of the authors (E.~Paal) would like to thank Jim Stasheff for reading
of the preliminary manuscript and useful discussions and suggestions.

\vskip1cm
\noindent
Department of Mathematics,
Tallinn Technical University\\
Ehitajate tee 5, 19086 Tallinn, Estonia\\
e-mails: liivi.kluge.ttu@mail.ee and epaal@edu.ttu.ee

\end{document}